# Buone pratiche didattiche per prevenire falsi positivi nelle diagnosi di discalculia: il progetto "PerContare"

# Preventing false positives in the diagnosis of dyscalculia through good teaching practices: the "PerContare" project

Anna E. Baccaglini-Frank[a], Maria G. Bartolini Bussi[b]

[a] Università degli Studi di Roma "La Sapienza", baccaglinifrank@mat.uniroma1.it
[b] Università degli Studi di Modena e Reggio Emilia, mariagiuseppina.bartolini@unimore.it

Abstract

Per contrastare il fenomeno dei falsi positivi nelle diagnosi di discalculia in Italia, tra i bambini di terza elementare, è stato realizzato un progetto triennale (2011-2014), fondato sulla collaborazione tra ricercatori in psicologia e didattica della matematica. Durante il progetto sono state elaborate e sperimentate pratiche didattiche atte a prevenire e affrontare difficoltà di apprendimento nell'ambito dell'aritmetica all'inizio della scuola primaria. In questo articolo vengono discusse le basi su cui si fonda il progetto, in particolare le basi scientifiche dei materiali didattici elaborati, di cui verranno proposti alcuni esempi. Tali esempi sono prototipici rispetto all'attenzione verso un approccio cinestetico-tattile e visuo-spaziale alla matematica. Le pratiche didattiche proposte hanno portato gli studenti ad interiorizzare relazioni parte-tutto e a pensare alla moltiplicazione in maniera strutturale, attraverso appositi diagrammi. Saranno forniti dati quantitativi e discussi brevemente dati qualitativi che confermano l'efficacia delle pratiche didattiche proposte.

Parole chiave: calcolo a mente; didattica inclusiva; difficoltà di apprendimento in matematica; discalculia; falsi positivi.

Abstract

To contrast the phenomenon of false positives in the diagnoses of dyscalculia in Italy, among 3rd grade children, a 3-year project (2011-2014), built upon a collaboration between psychologists and mathematics educators, was carried out. Within the project specific teaching strategies for preventing and addressing learning difficulties in arithmetic arising at the beginning of primary school were designed and tested in schools nation wide. This paper presents the project's background and the theoretical foundations of the didactical material designed, providing also some examples. In particular, prototypical examples will show how the activities, grounded within a kinesthetic-tactile and visual-spatial approach, are designed to lead to students' interiorization of part-whole relations, and to their thinking about multiplication through structured diagrams. Qualitative data confirming effectiveness of the proposed didactical strategies will also be discussed.

## 1. Che cos'è la discalculia: una questione delicata

La maggior parte della letteratura sulla discalculia nell'ambito della psicologia cognitiva studia lo sviluppo tipico di meccanismi di base per il processamento del numero, introducendo termini come "developmental dyscalculia", "mathematical learning disability (or disorder)" (l'acronimo più usato è MLD), "persistent mathematical learning difficulty" (in italiano, "discalculia", "discalculia evolutiva", "disturbo specifico del calcolo", "difficoltà persistenti in matematica"), tra altri per descrivere situazioni di sviluppo atipico (Passolunghi & Siegel, 2004; Piazza et al., 2010). Tuttavia le definizioni di tali costrutti sono ancora argomento di discussione, e la terminologia non è ancora utilizzata in modo coerente (Mazzocco, 2005; Mazzocco & Räsänen, 2013).

È possibile che tali incoerenze sussistano anche per il grande numero di ipotesi che sono state avanzate per spiegare prestazioni sotto la norma in matematica. Inoltre, negli studi sull'efficacia di interventi di recupero, in ambito psicologico, gli studenti sono soggetti a cui vengono somministrati esercizi, e le prove pre e post misurano soltanto particolari aspetti cognitivi, senza tener conto di altri fattori legati alla competenza matematica degli studenti. Tra i fattori generalmente esclusi, che sappiamo invece essere fondamentali per studiare il profilo di competenza matematica di uno studente, ci sono, per esempio, l'uso del linguaggio matematico, la partecipazione a discussioni matematiche, i significati costruiti durante l'interazione con ambienti digitali, e aspetti affettivi. Tutti questi elementi sono, invece, oggetto di ricerche specifiche nella didattica della matematica (Bartolini Bussi & Mariotti, 2009; Heyd-Metzuyanim, 2013; Mazzocco & Myers 2003; Santi & Baccaglini-Frank, 2015; Zan, 2007).

### 1.1. Il ruolo fondamentale dello sviluppo del "senso del numero"

Il costrutto "senso del numero" è attualmente usato con significati diversi in ambiti differenti di ricerca da parte di psicologi e studiosi di didattica della matematica, e non gode ancora di un'interpretazione univoca. Nemmeno nell'ambito più ristretto della didattica della matematica il costrutto è ben definito. Tuttavia, ciò su cui sembra esserci maggiore consenso è che alcuni aspetti del senso del numero sembrano avere importanti implicazioni per la didattica. Lo sviluppo del senso del numero è visto come essenziale per l'apprendimento dell'aritmetica (Griffin, Case, & Siegler, 1994) e alcuni ritengono che sia centrale per poter percepire e diventare consapevoli di aspetti strutturali del numero (Mulligan & Mitchelmore, 2013). Queste sono ritenute abilità fondamentali anche per promuovere il pensiero algebrico precoce (Resnick, Bill, Lesgold & Leer, 1991; Schmittau, 2011). Invece, un deficit di senso del numero è stato associato a basso rendimento in matematica, o a quello che in psicologia si chiama *low numeracy* (Butterworth & Laurillard, 2010; Halberda, Mazzocco & Feigenson, 2008), o discalculia (Dehaene & Wilson, 2007; Landerl, Bevan & Butterworth, 2004).

Molti studi provenienti dagli ambiti delle neuroscienze, della psicologia dello sviluppo, e della didattica della matematica suggeriscono che per migliorare l'apprendimento del numero sia importante imparare ad usare le dita in diversi modi: nel conteggio e nella rappresentazione dei numeri (Baccaglini-Frank & Maracci, 2015; Brissiaud, 1992); ma anche per favorire lo sviluppo delle abilità motorie fini e della gnosia digitale, che è la consapevolezza delle mani e delle dita nelle loro posizioni normali (Butterworth, 1999; Gracia-Bafalluy & Noël, 2008; Noël, 2005). Si ritiene, quindi, che l'uso sistematico delle dita possa avere un effetto positivo sullo sviluppo del senso del numero, già in età prescolare.

Joanne Mulligan e il suo gruppo di ricerca (Mulligan & Mitchelmore, 2013; Papic, Mulligan & Mitchelmore, 2011) non parlano di senso del numero ma di consapevolezza di *pattern and structure,* cioè regolarità e struttura, suggerendo che in bambini dai quattro ai nove anni tale consapevolezza è fondamentale per l'apprendimento dei concetti, delle competenze e delle abilità matematiche (Mulligan, 2011). Infatti, in diversi studi longitudinali, i dati raccolti da tali autori hanno mostrato chiaramente come studenti con basso rendimento scolastico in matematica in generale producono rappresentazioni male organizzate da un punto di vista spaziale: tendono a contare soltanto per uno e a non essere in grado di visualizzare relazioni parte-tutto.

Si può notare, infine, come negli ultimi anni sia aumentato molto l'interesse per materiali curricolari sviluppati per le prime classi di scuola primaria, o addirittura per la scuola dell'infanzia; ne è testimonianza, per esempio, lo Studio 23-esimo dell'International Commission on Mathematical Instruction (ICMI), organizzato a giugno 2015, attorno al tema "Primary Mathematics study on whole numbers"[1].

**1.2. La situazione italiana**

In Italia, gli studenti con disturbi specifici dell'apprendimento (DSA), fra cui la discalculia, sono stimati fra il 3% e il 5% della popolazione scolastica (MIUR, 2011a), negli ultimi anni, sembra esserci una tendenza all'aumento di certificazioni (MIUR, 2011b), che, in alcune regioni, aumenta di molti punti percentuali le stime sopra citate. Cerchiamo di capire meglio questo fenomeno, in particolare per quanto riguarda la discalculia. Dati dell'International Academy for Research on Learning Disabilities (IARLD) indicano che circa il 2,5% della popolazione dovrebbe presentare discalculia, eventualmente in comorbidità con altri disturbi specifici, come per esempio la dislessia, e che gli studenti con discalculia evolutiva pura dovrebbero essere meno dell'1% della popolazione (Cornoldi & Lucangeli, 2004). Tuttavia, uno studio sulla popolazione scolastica elementare italiana indica che più del 20% dei bambini alla fine della classe quinta risulta essere positiva ai test usati a livello nazionale per diagnosticare la discalculia (Lucangeli, 2005). Purtroppo non sono disponibili dati più recenti su questa fascia della popolazione italiana. In ogni modo, stando alle stime del MIUR, si deve concludere che moltissimi degli studenti che risultano positivi, e che purtroppo vengono etichettati come "discalculici", non ottengono questo risultato per vera inadeguatezza dei loro meccanismi cognitivi innati ma per una molteplicità di altri motivi. Dunque possono essere definiti *falsi positivi*.

In Italia il fenomeno di etichettamento di tutti gli studenti che presentano basse prestazioni in matematica ha creato e continua a creare molta tensione e confusione rispetto a che cosa sia la discalculia. Queste incertezze producono per lo più disagio tra insegnanti, studenti e famiglie, anche perché per legge gli studenti certificati DSA hanno diritto a particolari strumenti compensativi e misure dispensative (L. n. 170/2010). Nell'adattare tali strumenti e misure ai singoli studenti poco o nessun supporto è offerto agli insegnanti, che, secondo il sistema inclusivo italiano[2], si trovano in classe contemporaneamente studenti con bisogni educativi spesso molto diversi tra loro.

---

[1] Per ulteriori informazioni visitare http://www.mathunion.org/icmi/conferences/icmi-studies/ongoing-studies/icmi-study-23/.

Un'ipotesi emersa da questo contesto è che la proposta di strategie didattiche adeguate sin dall'inizio della scuola primaria (o anche prima) possa limitare drasticamente la comparsa di falsi positivi nelle diagnosi di discalculia evolutiva[3]. Tale ipotesi ha dato origine al progetto "PerContare".

## 2. Il materiale didattico elaborato nel progetto "PerContare"

Il progetto "PerContare" è un progetto inter-regionale italiano di tre anni (2011-2014) mirato a sviluppare strategie didattiche inclusive e materiale curricolare per aiutare insegnanti dei primi anni della scuola primaria (nelle classi prime e seconde) ad affrontare ed alleviare possibili difficoltà di apprendimento dei bambini, specialmente di quelli a rischio rispetto ad una diagnosi di discalculia. Il progetto nasce da una fruttuosa collaborazione tra ricercatori in didattica della matematica e in psicologia che si sono proposti i seguenti obiettivi:

- elaborare e testare strategie didattiche e guide didattiche per insegnanti di matematica delle classi prime e seconde che aiutino tutti gli studenti a sviluppare adeguate competenze numeriche;
- dare supporto agli insegnanti nell'imparare ad utilizzare i materiali e le strategie proposte in modo efficace;
- sviluppare strumenti per un'identificazione precoce delle difficoltà in aritmetica, da utilizzare nelle classi prime e seconde;
- sviluppare strumenti di potenziamento per bambini che continuano a manifestare difficoltà nello sviluppo di competenze numeriche di base.

In questo contributo si illustreranno le basi teoriche su cui sono costruite le pratiche didattiche proposte nelle guide per insegnanti, con alcuni esempi, e si illustreranno i risultati di uno studio longitudinale svolto nell'ambito del progetto, per mostrare l'efficacia dei materiali rispetto agli obiettivi posti.

### 2.1. Basi teoriche su cui sono fondate le pratiche didattiche proposte in "PerContare"

Il progetto "PerContare" prende in considerazione i principali contenuti dell'ambito aritmetico proposti nelle "Indicazioni per il curricolo della scuola dell'infanzia e del primo ciclo d'istruzione" (MIUR, 2012) per le classi prima e seconda. In tale contesto, i materiali fanno particolare attenzione a favorire lo sviluppo di consapevolezza di relazioni strutturali, in particolare la relazione parte-tutto (Baccaglini-Frank, 2015). Infatti, la consapevolezza di questo tipo di relazione è stata studiata da Resnick et al. (1991) in termini di schemi protoquantitivi che "organizzano le conoscenze dei bambini sui modi in cui il materiale che li circonda si ripartisce e si ricombina" (p. 32). Se l'attenzione viene portata su tale relazione, l'addizione e la sottrazione non sono più viste come operazioni separate, ma come azioni in un rapporto dialettico derivante proprio dalla relazione parte-tutto tra quantità (Schmittau, 2011). Inoltre, grazie alla consapevolezza di relazioni di tipo parte-tutto, gli studenti possono vedere i numeri come unità astratte che possono essere ripartite e ricombinate in modi diversi per facilitare il calcolo numerico. Un'attenzione alla relazione parte-tutto è altamente risonante con la teoria di Mulligan e colleghi (Mulligan & Mitchelmore, 2013; Papic et al., 2011) sul ruolo fondamentale della percezione di regolarità e struttura.

---

[2] Si vedano, per esempio, MIUR (2009), L. n. 170/2010, Direttiva Ministeriale 27 dicembre 2012.

[3] Ricordiamo che tale disturbo non può, per legge, essere diagnosticato prima della terza primaria o degli otto anni.

I riferimenti principali per le scelte metodologiche (Baccaglini-Frank & Bartolini Bussi, 2012; Baccaglini-Frank & Scorza, 2013a; 2013b) sono i seguenti: la teoria sui canali privilegiati per l'accesso e la produzione delle informazioni (Stella & Grandi, 2011), il modello del triplo codice di Dehaene (1992), e la teoria della mediazione semiotica (Bartolini Bussi & Mariotti, 2009). Si sottolinea, quindi, il riferimento sia a studi di psicologia che di ricerca in didattica della matematica, proprio per rispondere ai problemi di separazione tra le due aree di ricerca segnalate nel paragrafo 1.

La ricerca in psicologia cognitiva (Mariani, 1996; Stella & Grandi, 2011) ha identificato quattro canali specifici per l'accesso all'informazione e la produzione dell'informazione: il canale visivo verbale (linguaggio scritto), il canale visivo non-verbale o visuo-spaziale (immagini), il canale acustico e verbale (linguaggio parlato, suoni), e il canale cinestetico-tattile (manipolazione manuale). Per gli studenti con DSA risultano molto spesso essere privilegiati alcuni canali rispetto ad altri. Questi sono, in generale, i canali visivo non-verbale e cinestetico-tattile. Studi in didattica della matematica (Radford, 2014), confermati da studi sulle neuroscienze (Gallese & Lakoff, 2005), sottolineano che esperienze di tipo senso-motorio, percettivo, e cinestetico-tattili sono fondamentali nella formazione di concetti matematici, anche altamente astratti. Questi studi confermano, per esempio, il ruolo fondamentale di particolari modi di usare le dita per lo sviluppo del senso del numero (si veda il sottoparagrafo 1.1) e l'importanza della manipolazione di oggetti concreti che consentano di imparare facendo.

Rispetto all'uso di oggetti concreti, artefatti, nell'apprendimento della matematica la teoria della mediazione semiotica (Bartolini Bussi & Mariotti, 2009) getta luce su come l'uso di particolari artefatti, appositamente pensati e costruiti dalla società, da parte dello studente per risolvere problemi matematici possa favorire lo sviluppo di significati coerenti con i significati matematici che costituiscono l'obiettivo didattico di una certa lezione o di una serie di lezioni.

**2.2. Le attività nelle guide per insegnanti**

In questa sezione sono presentati alcuni esempi tratti dalle guide per gli insegnanti. La versione della guida a cui si fa riferimento è quella finale, che raccoglie in maniera ordinata sia i materiali proposti agli insegnanti nel corso della sperimentazione, sia la documentazione dell'attività svolta nelle classi coinvolte.

Ciascuna attività[4] è presentata secondo lo schema seguente: è data una stima della durata della lezione; si descrivono i materiali necessari per preparare la lezione, e sono riportate le consegne specifiche da dare agli studenti; se l'attività prevede diverse fasi di lavoro, queste sono descritte nei paragrafi successivi; segue una sezione che riporta ciò che l'insegnante si può aspettare, sulla base della sperimentazione triennale già svolta dei materiali didattici; poi compare una sezione in cui sono esplicitati i significati matematici di cui l'attività intende favorire la costruzione; seguono proposte su come aiutare gli studenti a costruire i significati matematici, in particolare quando sono previste discussioni matematiche, presentando anche video di discussioni matematiche; e sono poi fornite schede di lavoro per gli studenti, in vista di successive esercitazioni ed eventuali compiti per casa; spesso compaiono, infine, riferimenti ad alcuni dei software sviluppati all'interno del progetto che possano rafforzare abilità sviluppate dall'attività. Dalla struttura delle sezioni di cui è composta la presentazione di ciascuna attività emerge chiaramente il riferimento al quadro teorico della mediazione semiotica.

---

[4] Tutte le attività sono accessibili gratuitamente, previa iscrizione, al sito http://percontare.asphi.it.

**2.3. Primo esempio: la notazione posizionale decimale**

Oltre alle dita e alla linea dei numeri per lavorare con relazioni tra i numeri entro il 10 e poi entro il 20[5], uno degli artefatti fondamentali proposti per la prima classe sono le cannucce in fascetti di dieci. La rappresentazione dei numeri in fascetti e cannucce sparse consente di vedere i numeri in un formato analogico, ma contiene anche i semi della notazione simbolica in cifre arabe (i fascetti corrispondono alle decine e le cannucce sparse alle unità). In una delle prime attività nella guida per la prima classe che propone l'uso dei fascetti di cannucce si chiede al bambino di attuare passaggi di transcodifica da un codice all'altro del modello di Dehaene (1992): il codice verbale orale e scritto (per es. "tredici"), il codice simbolico scritto (per es. "13"), il codice analogico (per es. 13 stanghette). Questa attività è progettata per aiutare gli studenti a mettere in relazione i tre codici, rafforzando il legame in particolare tra i codici simbolico e analogico. I numeri in fascetti di cannucce mantengono, infatti, un forte legame concreto con una rappresentazione analogica accessibile e riproducibile usando il canale cinestetico-tattile, ma possono essere usati come trampolino per vedere e gestire la decomposizione in unità, decine (ed eventualmente centinaia, se si legano dieci fascetti) dei numeri naturali.

Le attività successive propongono l'uso di scatole trasparenti per decomporre e manipolare i numeri in unità e decine. A differenza dell'abaco, però, le scatole trasparenti consentono di lavorare con anche più di nove unità contemporaneamente, senza costringere il/la bambino/a ad interrompere il proprio flusso di pensieri e gestire il problematico passaggio da dieci unità ad una decina fino a quando è pronto/a. Inoltre, tale passaggio non richiede l'uso di astrazione nell'assumere che improvvisamente una pallina che prima valeva 1, ora vale 10. Infatti la decina di cannucce è costituita proprio dalle dieci cannucce sparse considerate, legate con un elastico, ed inserite in questo momento nella scatola delle decine (nella Figura 1 si vede una schermata del software Cannucce[6] che simula questo modello). Addizione e sottrazione possono essere presentate insieme, come nozioni complementari, che hanno senso soltanto se messe in relazione reciprocamente.

---

[5] Per un approfondimento si vedano Baccaglini-Frank (2015), Baccaglini-Frank e Scorza (2013b).

[6] Il software è stato ideato all'interno del progetto "PerContare" e implementato da Ivana Sacchi. È disponibile gratuitamente dal sito http://www.ivana.it/j/.

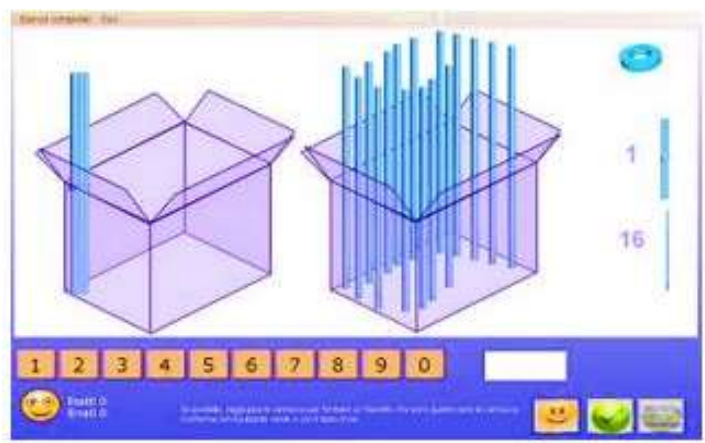

Figura 1. Schermata del software Cannucce sviluppato all'interno del progetto "PerContare" per consentire agli studenti di esercitarsi con il modello delle cannucce in fascetti e delle scatole trasparenti.

**2.4. Secondo esempio: la moltiplicazione con i diagrammi-rettangolo**

In questo caso si mostra un approccio alle cosiddette *tabelline*. Nella guida per insegnanti per la classe seconda si introduce l'artefatto *diagrammi-rettangolo* (Figura 2). Un diagramma rettangolo è un cartoncino (ed eventualmente una rappresentazione grafica su carta o sulla lavagna) costituito da una griglia di quadretti da 1 cm$^2$: il numero di righe rappresenta un fattore e il numero di colonne il secondo del prodotto rappresentato. Ogni diagramma-rettangolo (entro il $10 \times 10 = 100$) ha una collocazione all'interno della "casa dei rettangoli"[7] (Figura 3).

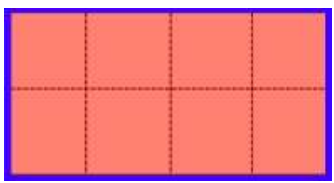

Figura 2. Il diagramma-rettangolo "$4 \times 2$" o "4 per 2 (volte)"[8].

Sono proposte diverse attività con i diagrammi-rettangolo con l'obiettivo didattico di favorire lo sviluppo di strategie visuo-spaziali e cinestetico-tattili per il calcolo di prodotti in base a fatti già conosciuti. In questo modo vengono ricavati tutti i prodotti entro il

---

[7] Tale tavola non si trova spesso nei libri di testo per la scuola primaria. Un'eccezione degna di nota è il libro "Matematica che passione!" (Bruscaglioni, Ferri & Mattiassich, 2000) in cui la tavola (rovesciata rispetto a quella di "PerContare") compare nella parte di "giochi e indovinelli 3".

[8] Una volta che gli studenti scoprono la simmetria lungo la diagonale della tavola che parte dall'angolo in basso a sinistra e arriva nell'angolo in alto a destra, che corrisponde alla proprietà commutativa della moltiplicazione, non è più importante insistere sulla distinzione tra righe e colonne, e il rettangolo "n×m" sarà identificabile anche come "m×n".

$10 \times 10=100$. Infatti, nella classe seconda gli studenti conoscono le sequenze di multipli del 2, del 5, del 10, e spesso quelli del 3, che possono quindi utilizzare per ricavare i prodotti mancanti mediante la scomposizione e la composizione.

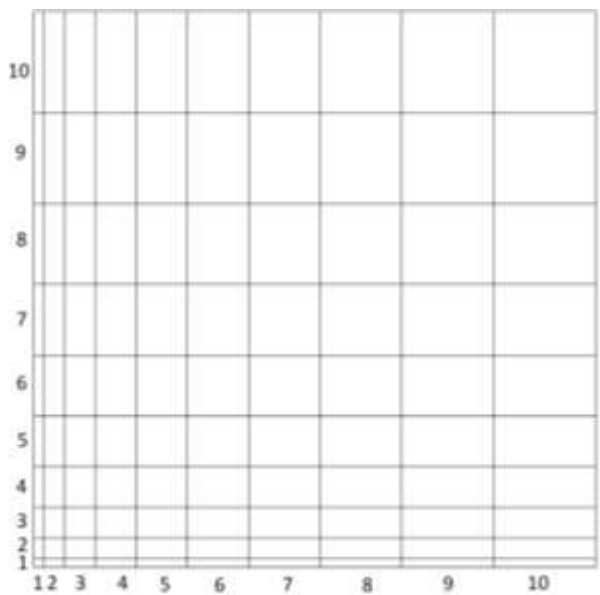

Figura 3. La "casa dei diagrammi-rettangolo".

Per esempio, per trovare il prodotto $7 \times 3$, i bambini possono utilizzare i diagrammi-rettangolo $5 \times 3$ e $2 \times 3$, avendo spezzato il 7 in 5 e 2 (Figura 4), e quindi sommando i prodotti conosciuti 15 e 6 mentalmente.

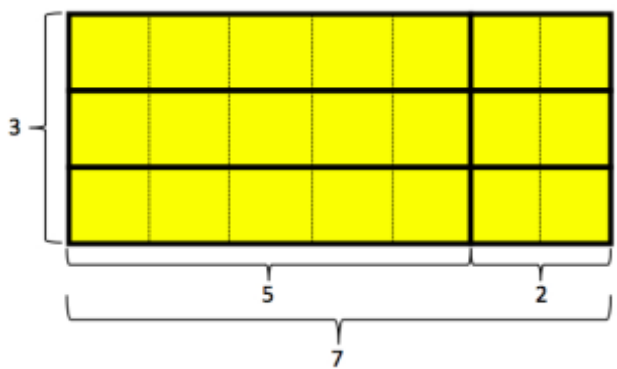

Figura 4. Una possibile scomposizione di $7 \times 3$ in $5 \times 3 + 2 \times 3$ con un diagramma-rettangolo.

## 3. Uno studio svolto nell'ambito del progetto

Durante il progetto è stato condotto uno studio longitudinale volto a valutare l'efficacia delle strategie didattiche implementate, rispetto all'obiettivo principale del progetto di prevenire i falsi positivi nelle diagnosi di discalculia. Lo studio si è svolto secondo la seguente metodologia.

È stato scelto un campione di 208 bambini (dieci classi) all'inizio della classe prima, che è stato seguito fino alla classe terza. Il campione non conteneva bambini con QI sotto norma ed era costituito da due gruppi: un gruppo sperimentale di 100 bambini (cinque classi) i cui insegnanti conoscevano ed utilizzavano i materiali di "PerContare", e un gruppo di controllo di 108 bambini (cinque classi) i cui insegnanti non conoscevano i materiali didattici di "PerContare". Le classi di controllo sono state scelte nello stesso ambiente socio-economico-culturale (stesso quartiere, talvolta una succursale della stessa scuola) con insegnanti non coinvolti in questa sperimentazione né in altre sperimentazioni sulla didattica della matematica. La scelta è stata operata dal gruppo di ricerca in collaborazione con i dirigenti scolastici. Si trattava di classi in cui era offerta una didattica di tipo tradizionale, senza la forte presenza degli artefatti scelti per le classi sperimentali. Le classi sperimentali hanno ricevuto in dotazione gli artefatti necessari alla sperimentazione (righelli, abaci, pascaline[9], bee-bot[10]) o indicazioni per costruirli (contamani, linee dei numeri[11], diagrammi-rettangolo). Gli insegnanti delle classi sperimentali parlavano ogni 15 giorni via Skype o di persona con la responsabile della guida, la prima autrice di questo lavoro, per discutere le loro idee rispetto all'implementazione dei materiali. Inoltre erano stati informati che i materiali erano già stati sperimentati l'anno prima in uno studio pilota e avevano accesso a materiale video girato nelle classi dell'anno precedente. In tutte le classi del campione sono state somministrate prove sulle abilità legate al numero e al calcolo, usati per la diagnosi di bambini a rischio (Biancardi, Mariani & Pieretti, 2011). Inoltre prove sviluppate nell'ambito di "PerContare" sono state somministrate nelle classi di entrambi i gruppi collettivamente sotto forma di gioco, a maggio in prima, a gennaio-febbraio e di nuovo a maggio in seconda.

Le prove sviluppate nell'ambito di "PerContare" per la classe prima valutano: (1) scrittura di numeri (entro il 1000, dettati non in ordine); (2) subitizing (con numerosità da 2 a 7); (3) stima (paragone di due numerosità); (4) enumerazione (conteggio di pallini e scrittura in notazione simbolica); (5) giudizio di grandezza (scegliere il simbolo che indica la quantità maggiore); (6) giudizio di quantità (decidere se due scritture, uno in formato analogico e l'altro in formato simbolico, rappresentano o meno lo stesso numero); (7) inserzioni sulla linea dei numeri (posizionare un numero su una linea dei numeri da 0 a 20 con tacche); (8) conteggio all'indietro (scrittura di numeri in ordine inverso sulle tacche di una linea dei numeri, a partire da un numero dato); (9) addizioni (operazioni scritte di cui tre necessitano di composizione di una decina); (10) sottrazioni (operazioni scritte in cui il numero maggiore è entro il 10). Per ciascuna consegna delle prove è stato calcolato il numero di risposte corrette.

Nella batteria per la classe seconda somministrata a gennaio-febbraio sette delle consegne sono dello stesso tipo rispetto alla batteria per la prima (le consegne 1, 2, 3, 5, 8, 9, 10), soltanto rese più complesse per la scelta dei numeri utilizzati, e altre tre consegne sono nuove (decomposizione, ordinamento crescente e decrescente). A maggio della classe seconda la prova utilizzata è come quella precedente per la stessa classe, con l'aggiunta di una consegna sulla moltiplicazione. Per ciascuna consegna delle prove è stato calcolato il numero di risposte corrette.

---

[9] Si veda Linea Scuola Quercetti al link
http://issuu.com/arcastudio/docs/cat_scuola2011_12?viewMode=magazine&mode=embed.

[10] Si veda http://www.campustore.it/bee-bot.

[11] Si vedano Baccaglini-Frank & Bartolini Bussi (2012), Baccaglini-Frank e Scorza (2013a; 2013b), Baccaglini-Frank (2015).

Per verificare la validità dei risultati ottenuti con le prove di "PerContare", a novembre della classe terza, al campione è stata somministrata anche la batteria AC-MT (Cornoldi, Lucangeli & Bellina, 2012).

## 4. Analisi dei dati raccolti

Lo studio ha previsto la raccolta di dati qualitativi (video, trascrizioni di episodi, schede di lavoro dei bambini) e di altri dati usati per un'analisi quantitativa (risultati delle prove proposte parallelamente nelle classi sperimentali e di controllo). Come testimonianza di come diversi bambini abbiano reagito alle attività proposte, riportiamo due episodi tra i dati qualitativi. Infine analizzeremo i dati quantitativi finora raccolti (alcuni sono ancora in fase di elaborazione).

### 4.1. La moltiplicazione con i diagrammi-rettangolo

I due episodi riguardano l'ambito del calcolo per mezzo di diagrammi-rettangolo e sono molto rappresentativi di un grande numero di episodi simili avvenuti nelle classi sperimentali, che, invece, erano molto più rari nelle classi del gruppo di controllo.

A marzo della classe seconda, durante la condivisione delle strategie utilizzate per calcolare $8 \times 6$ con i diagrammi-rettangolo, Chiara inventa il termine "quadretti fantasma"[12]. Dice che l'hanno aiutata "a semplificare il calcolo, ma poi bisogna ricordarsi di toglierli". Usa i "quadretti fantasma" per vedere l'8 come parte di 10, e ricavarsi il prodotto $10 \times 6$ (=60), da cui poi sottrae $2 \times 6$ (=12). Nonostante la proprietà distributiva non sia mai stata affrontata formalmente, Chiara, come molti altri bambini, è a suo agio nell'usarla per manipolare numeri.

Un secondo episodio avvenuto in un'altra seconda del gruppo sperimentale, ad aprile, è il seguente.

Insegnante: "Senza disegnare il palazzo[13] sette per otto, mi puoi dire come faresti a spezzarlo e a calcolarlo?"

Marco: "Allora, sette per otto…lo spezzo in cinque…e due e lo conto: cinque, dieci, quindici, venti, venticinque, trenta, trentacinque, quaranta [fa un gesto ripetitivo con cui alza in successione 8 dita]…e ho già quaranta. Poi conto i due: due, quattro, sei, otto, dieci, dodici, quattordici, sedici [ripete i gesti del conteggio come prima]. Poi faccio quaranta più, uh, spezzo il sedici in dieci e sei, e faccio quaranta più dieci cinquanta e poi aggiungo quei sei e ho cinquantasei."

Insegnante: "Caspita, ma sei bravissimo!"

Mentre parla, Marco più volte fissa il vuoto, come se stesse vedendo il palazzo che sta spezzando, contando, e ricomponendo. Procedure che fanno uso della visualizzazione

---

[12] Questo può anche essere visto come *segno pivot* secondo la Teoria della Mediazione Semiotica (Bartolini Bussi & Mariotti, 2009), ed è stato sfruttato come tale nell'intervento successivo dell'insegnante.

[13] In questa classe l'insegnante aveva accettato la proposta dei bambini di chiamare la tavola un insieme di "palazzi" in cui stavano i "Mr. Rettangolo" e questi si potevano "spezzare" e ricomporre.

come in questo esempio sono state rilevate molto frequentemente nelle classi sperimentali.

### 4.2. Analisi quantitativa dei dati delle prove somministrate nelle classi sperimentali e in quelle di controllo

I risultati delle prove somministrate alla fine della classe prima mostrano prestazioni paragonabili o superiori a quelle del gruppo di controllo in tutte le prove. Le prestazioni sostanzialmente superiori si sono verificate nelle seguenti consegne: giudizio di grandezza, addizione e sottrazione. Inoltre, nel gruppo sperimentale 4 dei 100 bambini hanno mostrato prestazioni scarse in almeno quattro delle consegne, mentre tale condizione si è verificata per 8 dei 108 studenti del gruppo di controllo. I risultati delle prove di gennaio-febbraio in seconda hanno confermato prestazioni significativamente superiori del gruppo sperimentale sulle addizioni e sottrazioni, e sulle consegne sull'ordinamento crescente e decrescente. La terza somministrazione delle prove a maggio della seconda ha confermato questi risultati.

Per quanto riguarda la validità delle prove elaborate all'interno del progetto, si ha una correlazione significativa ($p<0,05$) tra le nuove prove e la batteria standardizzata. In particolare, l'affidabilità è maggiore ($\alpha=0,8$) per le consegne che valutano la conoscenza dei numeri. Il paragone tra i punteggi medi ottenuti dai due gruppi sulle prove standardizzate di calcolo (AC-MT) mostrano una differenza significativa (t student = $p<0,05$) sulla velocità, sulle operazioni scritte e sulla conoscenza dei numeri. Il gruppo sperimentale ha riportato valori medi superiori su tutte le consegne delle prove standardizzate. Inoltre, nel gruppo sperimentale, la percentuale di studenti con prestazioni sotto il punteggio di cut off sull'AC-MT è circa la metà di quella corrispondente degli studenti del gruppo di controllo (7% e 13%, rispettivamente). Nessuno studente del gruppo sperimentale è risultato discalculico puro nelle diagnosi effettuate in classe terza.

### 5. Conclusioni

Riflessioni qualitative anche sulle prove somministrate hanno portato ad altre due considerazioni. In primo luogo, si è potuto notare come tutti i bambini del campione sperimentale abbiano risposto a tutti gli item proposti, a differenza di quanto è accaduto tra i bambini del gruppo di controllo, nel quale sono stati rilevati molti casi di omissioni. Inoltre i bambini del gruppo sperimentale hanno mostrato una grande varietà di strategie diverse per il calcolo mentale, insieme a maggiori abilità di controllo sulle soluzioni, a differenza dei bambini del gruppo di controllo che hanno usato strategie altamente standardizzate. In secondo luogo, analizzando le prove per scoprire i tempi di automatizzazione dei fatti additivi e moltiplicativi, si nota che per i fatti moltiplicativi, c'è stato un ritardo di circa tre mesi in media, rispetto agli studenti con prestazioni alte nel gruppo di controllo. Questo risultato era prevedibile in quanto "PerContare" non chiede ai bambini di memorizzare tabelline da un giorno all'altro come compiti per casa, ma propone un lavoro molto più articolato e profondo sulla manipolazione dei numeri.

Nel complesso, tali risultati suggeriscono che i materiali didattici di "PerContare" possano contribuire sostanzialmente a limitare il fenomeno dei falsi positivi nelle diagnosi di discalculia. Ma se si può ridurre il numero dei bambini positivi ai test per la discalculia in modo così importante soltanto con una didattica attenta, c'è da chiedersi che cosa ci dicano davvero le prove che si usano quotidianamente per la diagnosi, e, più in generale che cosa sia la discalculia. Queste sono questioni ancora aperte in diversi ambiti della ricerca. Poiché il raggiungimento di risposte certe è ancora lontano, in ambito didattico sembra che possa essere più fruttuoso smettere di cercare di scoprire *chi* siano i discalculici, per etichettarli e concentrare invece l'attenzione su *perché* alcuni studenti falliscano in certi ambiti (quali?) della matematica e ricercare che cosa è possibile fare per evitare tale fallimento.

## 6. Ringraziamenti

## Bibliografia

Baccaglini-Frank, A.E. (2015). Preventing low achievement in arithmetic through the didactical materials of the "PerContare" project. In X. Sun, B. Kaur & J. Novotná (eds.), *ICMI Study 23 Conference Proceedings* (pp. 169-176). Macau: University of Macau.

Baccaglini-Frank, A.E., & Bartolini Bussi, M.G. (2012). The "PerContare" project: proposed teaching strategies and some developed materials. In F. Dellai, I.C. Mammarella & A.M. Re (eds.), *International Academy for Research on Learning Disabilities 36th Annual Conference* (pp. 194-196). Trento: Erickson.

Baccaglini-Frank, A.E., & Scorza, M. (2013a). Preventing Learning Difficulties in early arithmetic: the "PerContare" project. In T. Ramiro-Sànchez & M.P. Bermùdez (eds.), *Libro de Actas I Congreso Internacional de Ciencias de la Educatiòn y des Desarrollo* (p. 341). Granada: Universidad de Granada.

Baccaglini-Frank, A.E., & Scorza, M. (2013b). Gestire gli studenti con DSA in classe: uso delle mani e della linea dei numeri nel progetto "PerContare". *Quaderni GRIMeD, 1*, 183–190.

Baccaglini-Frank, A.E., & Maracci, M. (2015). Multi-touch technology and preschoolers' development of number-sense. *Digital Experiences in Mathematics Education. 1*(1), 7–27.

Bartolini Bussi, M.G., & Mariotti, M.A., (2009). Mediazione semiotica nella didattica della matematica: artefatti e segni nella tradizione di Vygotskij. *L'Insegnamento della Matematica e delle Scienze Integrate, 32,* 269–294.

Biancardi, A., Mariani, E., & Pieretti, M. (2011). *La discalculia evolutiva. Dai modelli neuropsicologici alla riabilitazione*. Milano: Franco Angeli.

Brissiaud, R. (1992). A toll for number construction: finger symbol sets. In J. Bidaud, C. Meljac & J.P. Fischer (eds.), *Pathways to number. Children's developing numerical abilities*. New Jersey, NJ: Lawrence Erlbaum Associates.

Bruscaglioni, L., Ferri, F., & Mattiassich, M. (2000). *Matematica che passione. Per la Scuola elementare*. Vol.2. Torino: Il Capitello.

Butterworth, B. (1999). *What counts – how every brain is hardwired for math*. New York, NY: The Free Press.

Butterworth, B., & Laurillard, D. (2010). Low numeracy and dyscalculia: identification and intervention. *ZDM Mathematics Education*, *42*, 527–539.

"Cannucce" software. http://www.ivana.it/j/ (ver. 15.12.2015).

Cornoldi, C., & Lucangeli, D. (2004). Arithmetic education and learning disabilities in Italy. *Journal of Learning Disabilities*, *37*(1), 42–49.

Cornoldi, C., Lucangeli, D., & Bellina, M. (2012). *Test AC-MT 6-11 – Test di valutazione delle abilità di calcolo e soluzione di problemi*. Trento: Erickson.

Dehaene, S. (1992). Varieties of numerical abilities. *Cognition*, *44*, 1–42.

Dehaene, S., & Wilson, A. (2007). Number sense and developmental dyscalculia. In D. Coch, G. Dawson & K.W. Fischer (eds.), *Human behavior, learning, and the developing brain. Atypical development* (pp. 212-238). New York, NY: The Guilford Press.

Direttiva Ministeriale 27 dicembre 2012. *Strumenti di intervento per alunni con Bisogni Educativi Speciali e organizzazione territoriale per l'inclusione scolastica*.

Gallese, V., & Lakoff, G. (2005). The brain's concepts: the role of the sensory-motor system in conceptual knowledge. *Cognitive Neuropsychology*, *22*(3), 455–479.

Gracia-Bafalluy, M.G., & Noël, M.P. (2008). Does finger training increase young children's numerical performance?. *Cortex*, *44*, 368–375.

Griffin, S.A., Case, R., & Siegler, R.S. (1994). Rightstart: providing the central conceptual prerequisites for first formal learning of arithmetic to students at risk for school failure. In K. McGilly (ed.), *Classroom lessons:integrating cognitive theory and classroom practice* (pp. 24-49). Cambridge, MA: MIT Press.

Halberda, J., Mazzocco, M.M.M., & Feigenson, L. (2008). Individual differences in non-verbal number acuity correlate with maths achievement. *Nature*, *455*(7213), 665–668. http://doi.org/10.1038/nature07246 (ver. 15.12.2015).

Heyd-Metzuyanim, E. (2013). The co-construction of learning difficulties in mathematics-teacher-student interactions and their role in the development of a disabled mathematical identity. *Educational Studies in Mathematics*, *83*(3), 341–368.

ICMI. International Commission on Mathematical Instruction. http://www.mathunion.org/icmi/conferences/icmi-studies/ongoing-studies/icmi-study-23/ (ver. 15.12.2015).

Landerl, K., Bevan, A., & Butterworth, B. (2004). Developmental dyscalculia and basic numerical capacities: a study of 8-9-year-old students. *Cognition*, *93*, 99–125.

Legge 8 ottobre 2010, n. 170. *Nuove norme in materia di Disturbi Specifici di Apprendimento in ambito scolastico.* http://www.istruzione.it/esame_di_stato/Primo_Ciclo/normativa/allegati/legge170_10.pdf (ver. 15.12.2015).

Linea scuola Quercetti. http://issuu.com/arcastudio/docs/cat_scuola2011_12?viewMode=magazine&mode=embed (ver. 15.12.2015).

Lucangeli, D. (2005). *National survey on learning disabilities.* Roma: Italian Institute of Research on Infancy.

Mariani, L. (1996). *Investigating learning styles.* Reprinted from *Perspectives, Journal of TESOL-Italy, XXI*(2) and *XXII*(1).

Mazzocco, M.M.M. (2005). Challenges in identifying target skills for math disability screening and intervention. *Journal of Learning Disabilities*, *38*(4), 318–323.

Mazzocco, M.M.M., & Myers, G.F. (2003). Complexities in identifying and defining mathematics learning disability in the primary school years. *Annals of Dyslexia*, *53*, 218–253.

Mazzocco, M.M.M., & Räsänen, P. (2013). Contributions of longitudinal studies to evolving definitions and knowledge of developmental dyscalculia. *Trends in Neuroscience and Education*, *2*(2), 65–73.

MIUR. Ministero dell'Istruzione, dell'Università e della Ricerca (2009). *Linee guida per l'integrazione scolastica degli alunni con disabilità.*

MIUR. Ministero dell'Istruzione, dell'Università e della Ricerca (2011a). *Dislessia, Gelmini presenta misure a favore di studenti con disturbi specifici di apprendimento (DSA) per scuola e università.* http://hubmiur.pubblica.istruzione.it/web/ministero/cs200711 (ver. 15.12.2015).

MIUR. Ministero dell'Istruzione, dell'Università e della Ricerca (2011b). *Studenti con disturbi specifici dell'apprendimento. Rilevazioni integrative a.s. 2010-2011.* http://hubmiur.pubblica.istruzione.it/web/istruzione/prot5140_10 (ver. 15.12.2015).

MIUR. Ministero dell'Istruzione, dell'Università e della Ricerca (2012). Indicazioni nazionali per il curricolo della scuola dell'infanzia e del primo ciclo d'istruzione. *Annali della Pubblica Istruzione.* No. Speciale. http://www.annaliistruzione.it/var/ezflow_site/storage/original/application/55f6425315450eb079ff3e4da917750c.pdf (ver. 15.12.2015).

Mulligan, J. (2011). Towards understanding the origins of children's difficulties in mathematics learning. *Australian Journal of Learning Difficulties*, *16*(1), 19–39. http://doi.org/10.1080/19404158.2011.563476 (ver. 15.12.2015).

Mulligan, J.T., & Mitchelmore, M.C. (2013). Early awareness of mathematical pattern and structure. In L. English & J. Mulligan (eds.), *Reconceptualizing early mathematics learning* (pp. 29-46). Dordrecht: Springer Science-Business Media.

Noël, M.P. (2005). Finger gnosia: a predictor of numerical abilities in children? *Child Neuropsychology*, *11*, 1–18.

Papic, M.M., Mulligan, J.T., & Mitchelmore, M.C. (2011). Assessing the development of preschoolers' mathematical patterning. *Journal for Research in Mathematics Education*, *42*(3), 237–269.

Passolunghi, M.C., & Siegel, L.S. (2004). Working memory and access to numerical information in children with disability in mathematics. *Journal of Experimental Child Psychology*, *88*, 348–367.

"PerContare". http://percontare.asphi.it (ver. 15.12.2015).

Piazza, M., Facoetti, A., Trussardi, A.N., Berteletti, I., Conte, S., Lucangeli, D., … Zorzi, M. (2010). Developmental trajectory of number acuity reveals a severe impairment in developmental dyscalculia. *Cognition*, *116*(1), 33–41.

Radford, L. (2014). Towards an embodied, cultural, and material conception of mathematics cognition. *ZDM Mathematics Education*, *46*(3), 349–361.

Resnick, L.B., Bill, V.L., Lesgold, S.B., & Leer, N.M. (1991). Thinking in arithmetic class. In B. Means, C. Chelemer & M.S. Knapp (eds.), *Teaching advanced skills to at-risk students* (pp. 27-53). Menlo Park, CA: SRI international.

Schmittau, J. (2011). The role of theoretical analysis in developing algebraic thinking: a Vygotskian perspective. In J. Cai & E. Knuth (eds.), *Early algebraization: a global dialogue from multiple perspectives* (pp. 71-86). Berlin: Springer.

Santi, G., & Baccaglini-Frank, A.E. (2015). Possible forms of generalization we can expect from students experiencing mathematical learning difficulties. *PNA, Revista de Investigaciòn en Didàctica de la Matemàtica*, *9(3)*, 217-243.

Stella, G., & Grandi, L. (2011). *Come leggere la dislessia e i DSA*. Milano: Giunti Editore.

Zan, R. (2007). *Difficoltà in matematica – Osservare, interpretare, intervenire*. Milano: Springer Verlag.